\documentclass[12pt,a4paper]{article}
\usepackage{latexsym,amssymb}
\usepackage{t1enc}

\usepackage{amsmath}
\usepackage{amsthm}
\usepackage[cp1250]{inputenc}
\usepackage[T1]{fontenc}
\usepackage[english]{babel}
\usepackage[dvips]{graphicx}
\usepackage[a4paper, left=3.2cm, right=2.8cm, top=2.3cm, bottom=3.4cm, headsep=1.2cm]{geometry}

\newtheorem{theorem}{Theorem}[section]
\newtheorem{lemma}[theorem]{Lemma}
\newtheorem{proposition}[theorem]{Proposition}
\newtheorem{corollary}[theorem]{Corollary}
\newtheorem{remark}[theorem]{Remark}

\def\R{\mathbb{R}}
\def\C{\mathbb{C}}

\def\Z{\mathbb{Z}}

\def\d{\mathrm{\,d}}

\def\part{\partial}

\def\Ind{\mathrm{Ind}}
\def\Deg{\mathrm{Deg}}
\def\conv{\mathrm{conv}}

\begin{document}
\title{\textbf{{A Landesman-Lazer type result for periodic parabolic problems on $\R^N$ at resonance}}}
\author{ Aleksander \'{C}wiszewski, Renata {\L}ukasiak\\
\emph{Faculty of Mathematics and Computer Science} \\
\emph{Nicolaus Copernicus University}\\
\emph{ul. Chopina 12/18, 87-100 Toruń, Poland} } \maketitle

\begin{abstract}
We are concerned with $T$-periodic solutions of nonautonomous parabolic problem of the form $u_t = \Delta u + V(x) u + f(t,x,u)$,
$t >0$, $x \in \R^N$, with $V \in L^\infty (\R^N)+L^p(\R^N)$, $p \geq N$ and
$T$-periodic continuous perturbation $f:\R^N\times \R\to \R$. The so-called resonant
case is considered, i.e. when ${\cal N}:=\mathrm{Ker} (\Delta + V) \neq
\{0\}$ and $f$ is bounded. We derive a formula  for the fixed point index of the associated translation along trajectories operator in terms of the Brouwer topological degree of the time average mapping $\hat f: {\cal N}\to {\cal N}$ being the restriction of $f$ to ${\cal N}$. By use of the formula and continuation techniques we show that Landesman-Lazer type conditions imply the existence of $T$-periodic solutions.
\end{abstract}

\section{Introduction}
We are interested in the existence of $T$-periodic solutions of the
following nonlinear parabolic equation
\begin{equation}\label{23082013-0153}
\left\{   \begin{array}{l}
\displaystyle{\frac{\part u}{\part t}} (x,t) = \Delta u (x,t) + V (x) u (x,t) + f(t,x,u(x,t)), \ t >0, \, x\in\R^N,\\
u(\cdot, t) \in H^1 (\R^N), \, t \geq 0,
\end{array}  \right.
\end{equation}
where $\Delta$ is the Laplace operator (with respect to $x$), $V=V_0 - V_\infty$, $V_0 \in L^p(\R^N)$, $N \leq p < +\infty$, $V_\infty \in
L^\infty (\R^N)$ and $V_\infty \geq \bar v_\infty >0$ for some real number $\bar v_\infty>0$. The function $f: [0,+\infty) \times \R^N \times \R \to \R$ is assumed to be continuous, bounded, $T$-periodic with respect to time,
i.e.
\begin{equation}  \label{12032014-1537}
f(t,x,u) = f (t+T,x,u) \qquad \textnormal{for all } t \geq 0, x\in
\R^N, u \in \R,
\end{equation}
and satisfies the following conditions for all $t,s \in
[0,+\infty)$, $x \in \R^N$, $u,v \in \R$,
\begin{equation}\label{12032014-1513}
f(t,x,0) \leq M(x), \  \
\end{equation}
for some $M\in L^2(\R^N)$;
\begin{equation}\label{12032014-1514}
|f(t,x,u)-f(s,x,v)| \leq (\tilde K (x) +
K(x)|u|)|t-s|^{\theta} + L(x)|u-v|
\end{equation}
where $\theta\in (0,1)$, $\tilde K\in L^2(\R^N)$,
$K,L\in L^p (\R^N)$, $p\geq N$.
In this paper we consider the so-called resonant case, i.e. when the linear part $\Delta + V$ of the right-hand side of the equation has nontrivial kernel. Our main result reads as follows.
\begin{theorem}\label{15032014-2356}
Let ${\cal N}:=\mathrm{Ker} (\Delta + V) \neq \{0\}$, where $V$ is as above, and suppose that
$f: [0,+\infty) \times \R^N \times \R \to \R$ satisfies the conditions $(\ref{12032014-1537})$, $(\ref{12032014-1513})$, $(\ref{12032014-1514})$ and either
\begin{equation}\label{12032014-1611}
\int_0^T \bigg( \int_{\{\phi >0 \}}\check{f}_+(t,x)\phi(x)dx
+\int_{\{\phi <0 \}}\hat f_-(t,x)\phi(x)dx\bigg)dt>0
\end{equation}
for any $\phi \in {\cal N} \setminus \{0\}$, where $\check{f}_+(t,x):= \liminf_{s \to +\infty} f(t,x,s)$ and $\hat{f}_-(t,x):= \limsup_{s \to -\infty} f(t,x,s),$ or
\begin{equation}\label{12032014-1618}
\int_0^T \bigg( \int_{\{\phi >0 \}}\hat{f}^+(t,x)\phi(x)dx
+\int_{\{\phi <0 \}}\check f^-(t,x)\phi(x)dx\bigg)dt<0
\end{equation}
for any $\phi \in {\cal N} \setminus \{0\}$, where $\hat{f}^+(t,x):= \limsup_{s \to +\infty} f(t,x,s)$ and $\check{f}^-(t,x):= \liminf_{s \to -\infty} f(t,x,s)$ . Then {\em (\ref{23082013-0153})} admits a $T$-
periodic solution $u \in C([0,+\infty), H^2(\R^N)) \cap C^1([0,+\infty), L^2(\R^N))$.
\end{theorem}
\noindent Assumptions (\ref{12032014-1611}) and
(\ref{12032014-1618}) will be reffered to as Landesman-Lazer type
conditions, which had been widely used in the literature in the
context of evolutionary partial differential equations -- see e.g.
\cite{Fucik-Mawhin}, \cite{Brezis-Nirenberg}, as well as recent
papers \cite{Cwiszewski-JDIE}, \cite{Kokocki} and
\cite{Kokocki-2013}. The novelty of this paper may be viewed in the
fact that we study the problem on an unbounded domain, which entails
a few issues concerning compactness. This is a continuation of the
recent paper \cite{Cwiszewski-RL}, where we studied the periodic
parabolic problem without resonance.
\begin{remark} {\em
To indicate a class of functions satisfying the assumptions (\ref{12032014-1513}) and (\ref{12032014-1514}) consider $f:[0,+\infty) \times \R^N \times \R \to \R$ given by
$$
f(t,x,u) := U(x,t) + g(W(x,t)u),
$$
where continuous functions $U,W:\R^N \times [0,+\infty)\to \R$ are
such that $U(x, t)\leq U_0 (x)$ for all $x\in\R^N$ and $t>0$ with some bounded $U_0\in
L^2(\R^N)$,  $W(x,t) \leq L(x)$ for all $x\in\R^N$ and $t>0$ with some some
$L\in L^p (\R^N)$ and there is $\theta\in(0,1)$  such that, for all $t,s\geq 0$ and $x\in\R^N$,  $|U(t,x)-U(s,x)|\leq L_U (x)|t-s|^\theta$ and $|W(t,x)-W(s,x)|\leq L_W(x)|t-s|^\theta$ with $L_U\in L^2(\R^N)$ and  $L_W\in L^p(\R^N)$. Here $g:\R\to \R$ is a bounded
Lipschitz function such that $g(0)=0$. Then the
assumptions (\ref{12032014-1513}) and (\ref{12032014-1514}) are
satisfied. \hfill $\square$} \end{remark} \indent  Clearly, the
partial differential problem (\ref{23082013-0153}) can be
transformed into an abstract parabolic problem
\begin{equation} \label{15072014-1712}
\dot{u}(t) = -{\bold A} u(t)+ {\bold F} (t, u(t)), \ \ t \geq 0,
\end{equation}
by setting ${\bold A}: D({\bold A}) \to L^2(\R^N)$, with $D({\bold
A}):= H^2(\R^N)$, ${\bold A}u:= - (\Delta + V) u$, $u\in D({\bold
A})$, and ${\bold F}:[0,+\infty) \times H^1(\R^N) \to L^2(\R^N)$ is
given by $[{\bold F}(t,u)](x):=f(t,x,u(x))$, $x\in\R^N$, $t\geq 0$.
By the standard existence and uniqueness theory for evolution
equations we can properly define the translation operator ${\bold
\Phi}_T: H^1(\R^N)\to H^1(\R^N)$, ${\bold \Phi}_T (\bar u):= u(T)$,
$\bar u\in H^1 (\R^N)$, where $u:[0,+\infty)\to H^1(\R^N)$ is the
solution of (\ref{15072014-1712}) with the initial condition
$u(0)=\bar u$. In order to find $T$-periodic solutions of
(\ref{23082013-0153}) we shall look for fixed points of
${\bold \Phi}_T$ by use of local fixed point index theory.\\
\indent Motivated by \cite{Brezis-Nirenberg}, \cite{Cwiszewski-JDIE}
and \cite{Kokocki}, we prove a resonant version of \emph{averaging
principle}. Roughly speaking, it states that the topological
properties of our equation can be described in terms of the average
function $\bar {\bold F}: {\cal N} \to {\cal N}$ of ${\bold F}$,
restricted to the kernel of operator ${\bold A}$, given by
\begin{equation}
\bar{{\bold F}}(u):= \frac{1}{T}\int_0^T {\bold P} {\bold F}(t,u)\,
\d t \nonumber
\end{equation}
where ${\bold P}: L^2 (\R^N) \to {\cal N}$ is the orthogonal projection onto the finite dimensional space ${\cal N}$ (see Remark \ref{22082014-1040}). The main difficulty here comes from the lack
of the compactness of the translation operator.
In contrary to problems on bounded domains, the translation operator ${\bold \Phi}_T$ for this problem  is not completely continuous. Therefore we shall need to prove that the translation operator ${\bold \Phi}_T$ belongs to the class of ultimately compact operators, for which fixed point index theory is already known.
To this end we shall use the tail estimates technique that comes from Wang \cite{Wang}, who was interested in existence of the global attractor for the reaction-diffusion equation on $\R^N$, and  was also applied by Prizzi in \cite{Prizzi} to investigate the persistence of invariant sets in parabolic equations on unbounded domains.

We start with a parameterized family of problems
\begin{equation}\label{06052014-1913}
\dot{u}(t)= -{\bold A}u(t)+ \epsilon {\bold F}(t, u(t)), \ \ t>0,
\end{equation}
where $\epsilon \in [0,1]$ is a parameter, and let
${\bold \Phi}_T^{(\epsilon)}: H^1(\R^N) \to H^1(\R^N)$ be the translation along
trajectories operator for (\ref{06052014-1913}).
We prove that if $U \subset {\cal N}$ and $W \subset {\cal N}^\perp$ are open bounded sets such
that $0 \in W$ and $\bar {\bold F} (u) \neq 0$ for $u \in \part U$ then, for small $\epsilon\in (0,1]$,
\begin{equation}
\Ind ( {\bold \Phi}_T^{(\epsilon)}, U \oplus W) = (-1)^{m_-(\infty)}
\Deg_B (\bar {\bold F}, U), \label{11052014-0230}
\end{equation}
where $\Deg_B$ stands for Brouwer's topological degree and
$m_-(\infty)$ is the sum of the total multiplicities of the negative eigenvalues of $-(\Delta+V)$.
Here we exploit the spectral properties of the
operator $\Delta+V$. By use of the spectral theory, one may show that the essential spectrum of $-(\Delta+V)$ is contained in the interval $[\bar v_\infty, +\infty)$,
 which means that the set $\sigma_{ess}(-\Delta-V) \cap (-\infty,0)$ is finite and it consists of isolated eigenvalues of finite multiplicity. Thus the number
$m_- (\infty)$ is finite. The straightforward conclusion from
(\ref{11052014-0230}) is that the nontriviality of $\Deg_B (\bar
{\bold F}, U)$ gives the existence of $T$-periodic solutions of
(\ref{06052014-1913}). In the next step, by use of a continuation
argument, we show that, under some \emph{a priori bounds} condition,
the fixed point index of $\Phi_T$ (with respect to sufficiently
large balls) is equal to, up to a sign, $\Deg_B (\bar {\bold F},
U)$.
Finally, we show that the so-called Landesman-Lazer type conditions imply that the mentioned \emph{a priori bounds} estimates hold.\\
\indent The paper is organized as follows. In Section 2, we briefly
recall basic definitions from ultimately compact maps theory and
fixed point index theory for such maps and abstract parabolic
problems. Section 3 is devoted to the ultimate compactness property
of the translation along trajectories operator. In Section 4, we
derive an averaging and continuation principles. Finally, in Section
5, we prove a Landesman-Lazer type criterion for the existence of
periodic solutions of problem (\ref{23082013-0153}).

\section{Preliminaries}
\noindent {\bf Notation}. Throughout the paper we use the following
notational conveniences. If $(X,\|\cdot\|)$ is a normed space, then,
for $x_0\in X$ and $r>0$, we  put $B_X(x_0,r):=\{x\in X \mid
\|x-x_0\|<r\}$. If $Y \subseteq X$ is a subspace and $U \subset Y$,
then $\overline{U}^Y$ and $\part_Y U$ stand for the closure and the
boundary of $U$ in $Y$, respectively, and by $\part U$ and
$\overline{U}$ we denote the boundary and the closure of $U$ in $X$.
$\conv\, V$ and $\overline{\conv}^{X}\, V$ stand for the convex hull
and the closed (in $X$) convex hull of $V\subset X$, respectively.
By $(\cdot,\cdot)_X$ is denoted the inner product in $X$.
Finally, by $Y^\perp$ we denote the orthogonal complement of a subspace $Y$ of $X$ equipped with the inner product.\\

\noindent {\bf Measure of noncompactness}. Suppose that $\Omega$ is a
bounded subset of a Banach space $X$. Denote
$$
\beta(\Omega):= \inf \{r>0\ |\ \Omega\  \textnormal{can be covered
with a finite number of balls in}\  X \textnormal{of radius}\ r\}.$$
Then $\beta(\Omega)$ is finite and is called the Hausdorff measure
of noncompactness. It is easy to prove the following properties:
\begin{enumerate}
\item[(i)] $\beta(\Omega)=0$ if and only if $\Omega$ is relatively
compact;
\item[(ii)] $\beta(\conv\, \Omega)=
\beta(\Omega)=\beta(\overline{\Omega})$;
\item[(iii)] If $\Omega_1, \Omega_2 \subset X$ are bounded and such
that $\Omega_1 \subset \Omega_2$, then $\beta(\Omega_1)\leq
\beta(\Omega_2)$.
\end{enumerate}
More details concerning properties of the measure of noncompactness
can be found in \cite{Akhmerov-et-al} or \cite{Deimling}.\\

\noindent {\bf Ultimately compact maps and fixed point index}. A map
$\Phi: D \to X$, defined on a subset $D$ of a Banach space $X$ is
said to be {\em ultimately compact} if, for some $V\subset X$, the
equality $\overline{\conv}\,\Phi(V\cap D)=V$ implies that $V$ is
relatively compact. We shall say that an ultimately compact map
$\Phi:\overline{U} \to X$, defined on the closure of an open bounded
set $U\subset X$, is called {\em admissible} if $\Phi(u)\neq u$ for
all $u\in \part U$.  By an {\em admissible homotopy} between two
admissible maps $\Phi_0, \Phi_1:\overline{U}\to X$ we mean a
continuous map $\Psi:\overline{U} \times [0,1]\to X$ such that
$\Psi(\cdot, 0) = \Phi_0, \, \, \Psi(\cdot, 1)=\Phi_1$,
$\Psi(u,\mu)\neq u$ for all $u\in \part U$ and $\mu\in [0,1]$, and,
for any $V\subset X$, if $\Psi( (V\cap \overline{U}) \times [0,1]) =
V$, then $V$ is ultimately compact. Then $\Phi_0, \Phi_1$ are called
homotopic. It is worth mentioning that for the class of ultimately
compact maps a fixed point index can be considered. Its construction
can be found in \cite[1.6.3 and 3.5.6]{Akhmerov-et-al}. The fixed
point index for the discussed class of mappings posses
characteristic properties in the theory of compact operators. Below
we briefly recall these properties.


\begin{proposition}\label{23082013-0220} The fixed point index for the class of ultimately compact maps has the following properties.\\
{\em (i)   (existence)} If $\Ind(\Phi, U)\neq 0$, then there exists $u\in U$ such that $\Phi(u)=u$.\\
{\em (ii)  (additivity)} If $U_1, U_2 \subset U$ are open and $\Phi
(u)\neq u$ for all $u\in \overline{U\setminus (U_1\cup U_2)}$, then
$$\Ind(\Phi, U)=\Ind(\Phi, U_1) + \Ind(\Phi, U_2).$$
{\em (iii) (homotopy invariance)} If $\Phi_0, \Phi_1:\overline{U}\to
X$ are homotopic then
$$\Ind(\Phi_0, U)=\Ind (\Phi_1, U).$$
{\em (iv)  (normalization)} Let $u_0\in U$ and $\Phi_{u_0}:\overline
U\to X$ be defined by $\Phi_{u_0}(u)=u_0$ for all $u\in
\overline{U}$. Then $\Ind (\Phi_{u_0},U)$ is equal $0$ if
$u_0\not\in U$ and $1$ if $u_0\in U$.
\end{proposition}
\begin{remark} {\em  If $\Phi: \overline{U} \to X$ is a compact map then
$\Ind(\Phi, U)$ is an equal to the Leray-Schauder index $\Ind_{LS}
(\Phi, U)$ (see e.g. \cite{Granas}).}
\end{remark}

\noindent {\bf Evolution problems with perturbed sectorial operators}. Let $A:D(A)\to X$ be a sectorial operator such that for some $a>0$, $A+a I$ has its spectrum  in the half-plane $\{z\in \C \mid
\mathrm{Re}\, z >0 \}$. Let  $X^\alpha$, $\alpha>0$, be the fractional power space determined by $A+aI$. It is well-known that there exists $C_\alpha>0$ such that,
for all $u\in X$ and $t>0$,
$$
\| e^{-tA} u\|_{\alpha} \leq C_{\alpha} t^{-\alpha} e^{at} \|u\|
$$
where $\{ e^{-tA} \}_{t\geq 0}$ is the semigroup generated by $-A$.
Let $F:[0,\omega) \times X^\alpha \to X$, $\omega>0$ be such that,
for all $R>0$ there exist $L>0$ and $\theta \in (0,1)$, such that,
for all $t,s\in [0,\omega)$ and $u,v\in B(0,R)$,
$$
\|F(t,u)-F(s,v)\|\leq C(|t-s|^\theta + \|u-v\|_\alpha)
$$
and there exists $C>0$ such that, for all $t\in [0,\omega)$
and $u\in X^\alpha$,
$$
\|F(t,u)\|\leq C(1+\|u\|_\alpha).
$$
For $\bar u \in X^\alpha$, consider the equation
\begin{equation}\label{28012014-1055}
\left\{
\begin{array}{l}
\dot u(t) = - A u(t) + F(t,u(t)), \ t\in [0, \omega),\\
u(0) = \bar u.
\end{array}
\right.
\end{equation}
By a solution of $(\ref{28012014-1055})$ we understand
$$u\in C([0,\omega),X^\alpha)\cap C((0,\omega), D(A))\cap C^1((0,\omega), X)
$$
such that $(\ref{28012014-1055})$ holds. By classical results (see \cite{Henry} or \cite{Cholewa}), the problem $(\ref{28012014-1055})$
admits a unique global solution $u\in C([0,\omega),X^\alpha)\cap
C((0,+\omega), D(A))\cap C^1((0,\omega), X)$. Moreover, it is known that $u$ being solution of $(\ref{28012014-1055})$ satisfies the
following Duhamel formula
\begin{equation}
u(t)= e^{-tA} u(0) + \int_0^t e^{-(t-s)A}F(s,u(s))\, \d s,\qquad t\in [0,\omega).
\nonumber
\end{equation}

We shall use the following refinement of the continuity property (see \cite[Th. 3.1]{Cwiszewski-RL}).
\begin{theorem}\label{30082013-0308}
Assume that mappings $F_n : [0,\omega)\times X^\alpha\to X$, $n\geq 0$,
have the following properties
$$
\| F_n (t,u) \| \leq C(1+\| u \|_\alpha) \mbox{ for } t\in[0,\omega),\,
u\in X^\alpha, \, n\geq 0,
$$
for any $(t,x)\in [0,\omega)\times X^\alpha$ there exists a neighborhood
$U$ of $(t,x)$ in $[0,\omega)\times X^\alpha$ such that for all
$(t_1,u_1),(t_2,u_2)\in U$
$$
\|F_n (t_1,u_1)-F_n(t_2,u_2) \| \leq L(|t_1-t_2|^{\theta}
+\|u_1-u_2\|_{\alpha})
$$
for some $L > 0$ and $\theta\in (0,1)$ and, for each $u\in
X^\alpha$,
$$
\int_{0}^{t} F_n (s,u) \d s\to \int_{0}^{t} F_0(s,u) \d s  \ \
\mbox{ in } \ \ X  \mbox{ as } n\to +\infty
$$
uniformly with respect to $t$ from compact subsets  of $[0,\omega)$. Let
$u_n: [0,\omega) \to X^\alpha$, $n \geq 1$ be solution of
$(\ref{28012014-1055})$ with $F:=F_n$ and $\bar u := \bar u_n$. If
$u_n(0) \to u_0(0)$ in $X$, then $u_n (t)\to u_0 (t)$ in $X^\alpha$
uniformly with respect to $t$ from compact subsets of $(0,\omega)$, where
$u_0 :[0, \omega ) \to X^\alpha$ is a solution of
\begin{equation}
\left\{   \begin{array}{l}
\dot u(t) = - A u(t) + F_0(t,u(t)), \ t\in (0,\omega), \nonumber\\
u(0) = \bar u_0. \nonumber
\end{array}  \right.
\end{equation}
\end{theorem}

\section{Ultimate compactness property of translation along trajectories operator}
Let ${\bold A}_0: D({\bold A}_0) \to X$ be a linear operator in the space $X:=L^2(\R^N)$
given by
$$
{\bold A}_0 u:=-\sum_{i,j=1}^{N} a_{ij} \frac{\part^2 u}{\part x_j \part x_i} , \mbox{ for } u\in D({\bold A}_0):=H^2(\R^N),
$$
where $a_{ij}\in\R$, $i,j=1,\ldots, N$, are such that
$$
\sum_{i,j=1}^{N} a_{ij}\xi_i \xi_j \geq 0
\mbox{ for any } \xi=(\xi_1,\ldots,\xi_N) \in \R^N
$$
and $a_{ij}=a_{ji}$ for $i,j=1,\ldots,N$. Then $\bold A_0$ is a
self-adjoint, positive and sectorial operator in $L^2(\R^N)$. Define
${\bold V}_0: D({\bold V}_0) \to L^2(\R^N)$, $D({\bold V}_0):=
H^1(\R^N)$ by
$$
[{\bold V}_0 u] (x) := V_0(x) u(x), \ \ x\in\R^N,
$$
where $V_0 \in L^p(\R^N)$, $N \leq p < +\infty$ and let
${\bold V}_\infty : L^2(\R^N) \to L^2(\R^N)$ be given by
$$
[{\bold V}_\infty u](x):= V_\infty(x) u(x), \ \ x\in\R^N,
$$
where $V_\infty \in L^\infty(\R^N)$ and $V_\infty \geq \bar v_\infty > 0$ for some positive  $\bar v_\infty$.  Let  ${\bold A} := \bold A_0 - \bold V_0 + \bold V_\infty$ and  ${\bold
F} :[0,+\infty) \times H^1(\R^N)\to L^2(\R^N)$ be such that there
are $C > 0$, $K\in L^2(\R^N)$, $L\in L^p(\R^N)$ and $\theta \in
(0,1)$ such that, for any $t, s \geq 0$ and $ u , v \in H^1 (\R^N)
$,
\begin{equation}
\|\bold F(t,u) - \bold F(s, v)\|_{L^2} \leq C
(1+\|u\|_{H^1})|t-s|^\theta + C \|u-v\|_{H^1},
\label{23052014-0558}
\end{equation}
\begin{equation} \label{23052014-0559}
|{\bold F}(t,u)(x)| \leq L(x) |u(x)| + K (x)
(1+\|u \|_{H^1}) \mbox{ for a.a. } x\in \R^N.
\end{equation}
For any $\bar u \in H^1(\R^N)$ consider the following equation
\begin{equation}\label{11112013-1916}
\left\{   \begin{array}{l}
\dot u(t) = - {\bold A} u(t) + {\bold F} (t,u(t)),
\ t > 0,\\
u(0) = \bar u.
\end{array}  \right.
\end{equation}
Due to standard results in theory of abstract evolution equations (see \cite{Henry} or \cite{Cholewa})
the problem (\ref{11112013-1916}) admits a unique solution $u \in C ([0,+\infty),H^1(\R^N))$ $\cap C ((0,+\infty), H^2(\R^N))$ $\cap C^1 ((0,+\infty), L^2(\R^N))$. Inspired by \cite[Prop. 2.2]{Prizzi-FM} we have the following compactness result (being a version of \cite[Lem. 4.3]{Cwiszewski-RL}).
\begin{lemma} \label{08052014-1217}
Let $T>0$ and suppose that $u(\cdot; \bar u):[0,T]\to H^1(\R^N)$ is the solution of {\em (\ref{11112013-1916})} such that $\|u(t;\bar u)\|_{H^1} \leq R$ for all
$t \in [0,T]$ and some fixed $R>0$. Then there exists
a sequence $(\alpha_n)$ with $\alpha_n \to 0$ as $n \to \infty$ such that
\begin{equation}
\int_{\R^N \setminus B(0,n)} |u(t;\bar u)|^2\, \d x \leq R^2 e^{-2 \bar v_\infty t}+\alpha_n \ \ \textnormal{ for all\ } \ \ t\in [0,T], n\geq 1,
\nonumber
\end{equation}
where $\alpha_n\geq 0$ depend only on $p$, $N$, $K$, $L$, $V_0$,
$V_\infty$ and $a_{ij}'s$.
\end{lemma}
\begin{remark}
{\em (a) If $f:[0,+\infty)\times \R\to \R$ satisfies (\ref{12032014-1513}) and (\ref{12032014-1514}), then one can directly verify that the Nemytskii operator ${\bold F}$
determined by $f$ (i.e. given by the formula $[{\bold F}(t, u)](x):= f(t,x,u(x))$) satisfies (\ref{23052014-0558}) and (\ref{23052014-0559})
(see \cite[Lem. 4.1]{Cwiszewski-RL}). Here, keeping in mind our further needs, we do not assume that ${\bold F}$ is a Nemytskii operator.\\
\indent (b) Clearly, condition (\ref{23052014-0559}) implies that
${\bold F}$ has a sublinear growth. \hfill $\square$ }\end{remark}

\noindent {\bf Proof of Lemma \ref{08052014-1217}}.
Denote $u:=u(\cdot, \bar u)$. Let $\phi:[0,+\infty)\to \R$ be a smooth function such that $\phi(s)\in [0,1]$ for $s\in [0,+\infty)$, $\phi_{|[0,1]}\equiv 0$ and $\phi_{|[2,+\infty)}\equiv 1$ and let $\phi_n:\R^N\to\R$ be defined by $\phi_n (x):=\phi(|x|^2/n^2)$,
$x\in\R^N$, $n\in \Z$, $n\geq 1$. Then, using the regularity of solution, for any $t \in (0,T]$ one has
\begin{eqnarray*}
\frac{1}{2}\frac{\d}{\d t} (u(t), \phi_n u(t))_{L^2} & =&
\frac{1}{2}\bigg( (u(t),\phi_n \dot u(t))_{L^2} + (\dot u(t), \phi_n u(t))_{L^2}\bigg) =  (\phi_n u(t),\dot u(t))_{L^2} \\
& = & I_1(t)+I_2(t) + I_3(t)
\end{eqnarray*}
where
\begin{eqnarray*}
& I_1(t)& := -( \phi_n u(t), {\bold A}_0 u(t) )_{L^2},\\
& I_2 (t)& := ( \phi_n u(t),  ({\bold V}_0 -{\bold V}_\infty) u(t))_{L^2},\\
& I_3(t)& := ( \phi_n u(t),{\bold F} (t,u(t)))_{L^2}.
\end{eqnarray*}
First observe that
\begin{align}
I_1(t) & = -\int_{\R^N}\sum_{i,j=1}^N a_{ij} \frac{\part}{\part
x_j}(\phi_n(x) u(t))\frac{\part}{\part x_i} (u(t)) \d x
\nonumber \\
& = -\int_{\R^N} \phi_n(x) \sum_{i,j=1}^N a_{ij}\frac{\part}{\part
x_j}(u(t))\frac{\part}{\part x_i}(u(t)) \d x \nonumber \\
& - \frac{2}{n^2} \int_{\R^N}\phi'(|x|^2/n^2) \sum_{i,j=1}^N a_{ij}
x_j \frac{\part}{\part x_i}(u(t)) u(t) \d x \nonumber \\
& \leq \frac{2 L_\phi}{n^2} \int_{\{n \leq |x| \leq
\sqrt{2}n\}}\sum_{i,j=1}^N a_{ij} |x| |u(t)| \nabla u(t)| \, \d x \nonumber \\
&\leq \frac{2\sqrt{2}  L_\phi M
N^2}{n}\|u(t)\|_{L^2}\|u(t)\|_{H^1},\label{20042014-1122}
\end{align}
where $L_\phi:= \sup_{s\in[0,+\infty)}  |\phi'(s)|< \infty$ (as
$\phi$ is smooth and $\phi'$ is nonzero on a bounded interval) and $M:= \max_{1
\leq i,j \leq N}|a_{ij}|$. To estimate the second term, we see that
\begin{align}
I_2 (t) & \leq - \bar v_\infty \int_{\R^N}\phi_n(x)|u(t)|^2 \d x + C^{N/p}
\|u(t)\|_{H^1}^2 \bigg(\int_{\{|x| \geq n\}} |V_0(x)|^p
\d x\bigg)^{1/p}, \label{03092014-2341}
\end{align}
where $C=C(N)>0$ is the constant related to the Sobolev embedding
$H^1(\R^N) \subset L^{2N/(N-2)} (\R^N)$. Finally ,
\begin{align}
I_3(t) & = \int_{\R^N} \phi_n(x)\, {\bold F} (t,u(t))(x) u(t)\, \d x \nonumber\\
 & \leq \int_{\{|x| \geq n\}} L (x) |u(t)|^2 \, \d x
+ (1+R) \int_{\{|x| \geq n\}} K (x)|u(t)|\, \d x \nonumber \\
& \leq R^2C^{N/p}\bigg(\int_{\{|x| \geq n\}} |L(x)|^p \d x\bigg)^{1/p}  +  R
(1+R)\bigg(\int_{\{|x| \geq n\}} |K(x)|^2 dx\bigg)^{1/2}
\label{20042014-1139}
\end{align}
where $C$ is as above. Combining (\ref{20042014-1122}),
(\ref{03092014-2341}) and (\ref{20042014-1139}), we get
\begin{equation}
\frac{d}{dt}\int_{\R^N} \phi_n(x) |u(t)|^2 \d x \leq
-2 \bar v_\infty  \int_{\R^N} \phi_n(x) |u(t)|^2 \d x + \alpha_n \nonumber
\end{equation}
where $(\alpha_n)_{n \in \mathbb{N}}$ is a sequence such that
$\alpha_n \to 0$ as $n \to +\infty$. Multiplying by $e^{2 \bar v_\infty t}$ and integrating over $[0, \tau]$ we have
\begin{equation}
\int_{\R^N}\phi_n(x)|u(t)|^2 dx \leq e^{-2 \bar v_\infty
t} \int_{\R^N}\phi_n(x)|u(0)|^2 dx + \alpha_n
\end{equation}
where $(2\bar v_\infty)^{-1}(1-e^{-2 \bar v_\infty t})\alpha_n$ is again denoted by $\alpha_n$. This finishes the proof. \hfill $\square$\\

\indent Now we are going to show that the translation operator is ultimately compact.  We shall consider a parameterized problem. Suppose that $a_{ij}\in C([0,1],\R)$, $i,j=1,\ldots, N$, are such that the ellipticity condition holds
$\sum_{i,j=1}^{N} a_{ij}(\mu)\xi_i \xi_j \geq 0$ for any $\xi\in \R^N$ and $\mu\in [0,1]$. Let ${\bold A}_{0}^{(\mu)}: D({\bold A}_{0}^{(\mu)}) \to L^2(\R^N)$ be given by $D({\bold A}_{0}^{(\mu)})=H^2(\R^N)$,
 $$ {\bold A}_{0}^{(\mu)} u:= - \sum_{i,j=1}^{N} a_{ij}(\mu) \frac{\part^2 u}{\part x_j \part x_i}.$$
Assume that ${\bold F}: [0,+\infty) \times H^1(\R^N) \times [0,1] \to L^2(\R^N)$
is such that there are $C > 0$, ${\bold B}\in {\cal L}(H^1(\R^N), H^1(\R^N))$, $K\in L^2(\R^N)$ and $L\in L^p(\R^N)$, $p\geq N$,
$\theta \in (0,1)$ such that, for any $t, s
\geq 0$, $ u , v \in H^1 (\R^N)$ and $\mu, \nu \in [0,1]$,
\begin{equation}
\|{\bold F}(t,u, \mu) - {\bold F} (s, v, \mu)\|_{L^2}
\leq C (1+\|u\|_{H^1})|t-s|^\theta + C \|u-v\|_{H^1};
\label{18042014-0941}
\end{equation}
\begin{equation}
 | {\bold F} (t,u,\mu) (x) | \leq L(x) | [{\bold B} u](x) | + K(x)
 (1+ \|u\|_{H^1});
\label{18042014-0949}
\end{equation}
\begin{equation}
\left\| {\bold F} (t, u, \mu)\!-\! {\bold F}
(t, u, \nu)\right\|_{L^2} \leq
\left|\rho(\mu)\!-\!\rho(\nu)\right| \left(1+\|u\|_{H^1}\right) \label{18042014-0950}
\end{equation}
where $\rho \in C([0,1], \R)$.  Consider the problem
\begin{equation}\label{29082013-1258}
\dot u (t) = - {\bold A}^{(\mu)} u(t)
+ {\bold F} (t,u(t),\mu), \ t>0,
\end{equation}
where ${\bold A}^{(\mu)}
:= {\bold A}_{0}^{(\mu)}-{\bold V}_0+{\bold V}_\infty$.
As before, due to classical results in theory of abstract evolution problems, the problem (\ref{29082013-1258}) admits a unique global solution
$u \in C([0,+\infty), H^1(\R^N)) \cap C((0,+\infty), H^2(\R^N)) \cap C^1((0,+\infty), L^2(\R^N))$. Denote
by $u(\cdot; \bar u, \mu)$ the solution of (\ref{29082013-1258}) satisfying the initial value condition $u ( 0 ) = \bar u$.  Slightly modifying the proof of Lemma 4.4 in \cite{Cwiszewski-RL} we obtain
the following result.
\begin{lemma}\label{26082013-2147}
Take any $\bar u_1, \bar u_2\in H^{1}(\R^N)$ and $\mu_1, \mu_2\in [0,1]$ and suppose that there are solutions $u(\cdot;\bar u_i,\mu_i):[0,T]\to H^1(\R^N)$, $i=1,2$ of {\em (\ref{29082013-1258})}, for some fixed $T>0$. If $\|u(t;\bar u_1,\mu_1)\|_{H^1}\leq R$ and $\|u(t;\bar u_2,\mu_2)\|_{H^1} \leq R$ for all  $t\in [0,T]$ and
some fixed $R>0$, then there exists a sequence $(\alpha_n)$ with
$\alpha_n\to 0$ as $n\to\infty$ such that
$$
\int_{\R^N \setminus B(0,n)} \left|u(t;\bar u_1,\mu_1)-u(t;\bar
u_2,\mu_2)\right|^2 \d x  \leq e^{-2 \bar v_\infty t} \| \bar u_1 -
\bar u_2 \|_{L^2}^{2}  + Q \eta(\mu_1,\mu_2) + \alpha_n,
$$
for all $t\in [0,T]$ and $n\geq 1$, where $\alpha_n\geq 0$ and $Q>0$
depend only on $p$, $N$, $K$, $L$, $V_0$, $V_\infty$, ${\bold B}$ and $a_{ij}'s$,
$$
\eta(\mu_1,\mu_2):= \max\left\{
|\rho(\mu_1)-\rho(\mu_2)|,\max_{i,j=1,\ldots, N}\{
|a_{ij}(\mu_1)-a_{ij}(\mu_2)|\}\right\}.
$$
\end{lemma}
\noindent {\bf Proof:}  Let $\phi_n:\R^N\to\R$, $n\geq 1$, be as in the proof of Lemma \ref{08052014-1217}. Put $u_1:=u(\cdot; \bar u_1,\mu_1)$, $u_2:=u(\cdot; \bar u_2,\mu_2)$ and $w:=u_1-u_2$.
Observe that
\begin{eqnarray*}
\frac{1}{2}\frac{\d}{\d t} (w(t), \phi_n w(t))_{L^2} & =&
\frac{1}{2}\bigg(
 (w(t),\phi_n \dot w(t))_{L^2} + (\dot w(t), \phi_n w(t))_{L^2}\bigg) =  (\phi_n w(t),\dot w(t))_{L^2} \\
& = & I_1(t)+I_2(t) + I_3(t) + I_4 (t)
\end{eqnarray*}
where
\begin{eqnarray*}
& I_1(t)& := \left( \phi_n w(t), - {\bold A}_{0}^{(\mu_1)}u_1(t) + {\bold A}_{0}^{(\mu_1)}u_2(t) \right)_{L^2},\\
& I_2 (t)& := \left( \phi_n w(t), -{\bold A}_{0}^{(\mu_1)}u_2(t) + {\bold A}_{0}^{(\mu_2)}u_2(t) \right)_{L^2},\\
& I_3(t)& := \left( \phi_n w(t), ({\bold V}_0 -{\bold V}_{\infty})w(t) \right)_{L^2},\\
& I_4(t)& := \left( \phi_n w(t),{\bold F} (t,u_1(t),\mu_1)- {\bold F} (t,u_2(t),\mu_2) \right)_{L^2}.
\end{eqnarray*}
As for the first term we notice that
\begin{eqnarray*}
I_1(t)& =&
\left( \phi_n w(t), - {\bold A}_{0}^{(\mu_1)} w(t) \right)_{L^2} \\
&   = &  -\int_{\R^N}\sum_{i,j=1}^{N} a_{ij}(\mu_1) \frac{\part
}{\part x_j}(\phi_n (x) w(t))  \frac{\part}{\part x_i}(w(t)) \d x\\
& = & - \int_{\R^N} \phi_n (x)
\sum_{i,j=1}^{N} a_{ij}(\mu_1) \frac{\part }{\part x_j}(w(t))
\frac{\part}{\part x_i}(w(t)) \d x \\
& & -\frac{2}{n^2} \int_{\R^N} \sum_{i,j=1}^{N} \phi'(|x|^2/n^2)w(t) x_j a_{ij}(\mu_1) \frac{\part}{\part x_i} (w(t)) \d x\\
& \leq & \frac{2 L_\phi }{n^2} \int_{\{n\leq |x|\leq \sqrt{2}n\}} \sum_{i,j=1}^{N} a_{ij}(\mu_1) |x| |w(t)||\nabla w(t)|\d x\\
& \leq & \frac{2\sqrt{2} L_\phi M N^2}{n}
\|w(t)\|_{L^2}\|w(t)\|_{H^1}
\end{eqnarray*}
where $M:=\max_{1\leq i,j \leq N, \, \mu\in [0,1]} |a_{ij}(\mu)|$
and $L_\phi:=\sup_{s \in [0,+\infty)}\left|\phi'(s)\right| <
+\infty$. Further, in a similar manner
\begin{eqnarray*}
I_2(t) & = & - \int_{\R^N}\sum_{i,j=1}^{N}  \frac{\part }{ \part x_j} (\phi_n w(t)) (a_{ij}(\mu_1)-a_{ij}(\mu_2)) \frac{\part }{\part x_i }(u_2(t)) \d x  \\
& = & - \int_{\R^N} \phi_n (x) \sum_{i,j=1}^{N} (a_{ij}(\mu_1)- a_{ij}(\mu_2)) \frac{\part }{\part x_j}( w (t) ) \frac{\part}{\part x_i}( u_2(t) ) \d x \\
& & -\frac{2}{n^2} \int_{\R^N} \sum_{i,j=1}^{N} \phi'(|x|^2/n^2) w(t) x_j (a_{ij}(\mu_1)- a_{ij}(\mu_2))  \frac{\part}{\part x_i} (u_2 (t)) \d x \\
& \leq & N^2 \eta(\mu_1,\mu_2) \|w(t)\|_{H^1} \|u_2(t)\|_{H^1} + \frac{4\sqrt{2} L_\phi  \eta(\mu_1,\mu_2) N^2}{n} \|w(t)\|_{L^2} \|u_2(t)\|_{H^1}
\end{eqnarray*}
and
\begin{align}
I_3 (t) & \leq - \bar v_\infty \int_{\R^N}\phi_n(x)|w(t)|^2 \d x + C^{N/p}
\|w(t)\|_{H^1}^2 \bigg(\int_{\{|x| \geq n\}} |V_0(x)|^p
\d x\bigg)^{1/p} \label{20042014-1130}
\end{align}
where $C$ is a constant of the embedding $H^1(\R^N)\subset L^{\frac{2N}{N-2}}(\R^N)$.
To estimate $I_4 (t)$ observe that
\begin{eqnarray*}
I_4 (t)& =  &  \int_{\R^N} \phi_n(x) \left( {\bold F}(t,u_1(t),\mu_1)\!-\! {\bold F}(t, u_2(t),\mu_1)\right) w(t) \d x\\
& & + \int_{\R^N} \phi_n(x) \left( {\bold F} (t,u_2(t),\mu_1) \!-\! {\bold F} (t, u_2(t),\mu_2)\right) w(t) \d x. \label{07052014-1613}
\end{eqnarray*}
Further, by the H\"{o}lder inequality, it follows that
\begin{align}
\ & \int_{\R^N} \phi_n(x) \left( {\bold F}(t,u_1(t),\mu_1)\!-\! {\bold F}(t, u_2(t),\mu_1)\right) w(t) \d x \nonumber \\
\ &  \leq \int_{\{|x|\geq n\}}
(\left|{\bold F} (t,u_1(t),\mu_1)|+|{\bold F} (t,u_2(t),\mu_1)| \right)\, |w(t)|\, \d x \nonumber \\
\ & \leq \int_{\{|x| \geq n\}}\!\!\! L (x)(|{\bold B}u_1(t)|+|{\bold B}u_2(t)|) |w(t)| \, \d x
+ 2(1+R)\int_{\{|x| \geq n\}} \!\!\!\! K (x) |w(t)| \, \d x \nonumber \\
& \leq  8R^2 C^{N/p}\|{\bold B}\|_{\mathcal{L}{\cal}(H^1,H^1)}^{N/p} \bigg(\int_{\{|x| \geq n\}} |L(x)|^p
\d x\bigg)^{1/p}  +  4R (1+R)\bigg(\int_{\{|x| \geq n\}} |K(x)|^2 \d x \bigg)^{1/2}
\label{07052014-1601}
\end{align}
where $C$ is as above. Using (\ref{18042014-0950}),
we have
\begin{equation}
\int_{\R^N}\!\! \phi_n(x) \left(
{\bold F}(t,u_2(t),\mu_1)\!-\! \bold F(t,u_2(t),\mu_2)\right) w(t) \d x \leq \left| \rho(\mu_1)\! -\! \rho(\mu_2) \right| \!(1\!+\!R)2R. \label{07052014-1609}
\end{equation}
Combining (\ref{07052014-1613}) and
(\ref{07052014-1609}) we obtain
\begin{align}
I_4(t) & \leq
\left| \rho(\mu_1)\! -\! \rho(\mu_2) \right| \!(1\!+\!R)2R +
4RC^{N/p}\bigg(\int_{\{|x| \geq n\}} |L(x)|^p \d x\bigg)^{1/p} \|{\bold B}\|_{\mathcal{L}{\cal}(H^1,H^1)}^{N/p}  \nonumber \\ & \ \ \ \  \ \ \ \ \  \ \ \ \  +  4R (1+R)\bigg(\int_{\{|x| \geq n\}} |K(x)|^2 dx\bigg)^{1/2}.
\nonumber
\end{align}
Summing up all the estimates, we get, for any $n\geq 1$,
$$
\frac{\d}{\d t} (w(t), \phi_n w(t))_{L^2} \leq - 2 \bar v_\infty (w(t), \phi_n w(t))_{L^2} + \tilde C\eta(\mu_1,\mu_2) + \alpha_n
$$
for some constant $\tilde C=\tilde C(p,N, V_0, V_\infty, K, L, R)>0$. Multiplying by $e^{2 \bar v_\infty t}$ and integrating over $[0,\tau]$ one obtains
$$
e^{2 \bar v_\infty \tau}(w(\tau), \phi_n w(\tau))_{L^2}- (w(0), \phi_n
w(0))_{L^2} \leq  (2 \bar v_\infty)^{-1}(e^{2 \bar v_\infty \tau}-1)\,
(\tilde C \eta (\mu_1,\mu_2) +\alpha_n),
$$
which gives
$$
(w(\tau), \phi_n w(\tau))_{L^2} \leq e^{-2 \bar v_\infty
\tau}\|w(0)\|^2_{L^2} + (2 \bar v_\infty)^{-1} \left( \tilde C \eta
(\mu_1,\mu_2)+\alpha_n\right).
$$
And this finally implies the assertion as
$\|\phi_n w(\tau)\|_{L^2}^{2}\leq (w(\tau), \phi_n w(\tau))_{L^2}$. \hfill $\square$\\

Let ${\bold \Psi}_t: H^1(\R^N)\times [0,1]\to H^1(\R^N)$, $t>0$, be  the
translation operator for (\ref{29082013-1258}), i.e.
$ {\bold \Psi}_t (\bar u,\mu) = u (t ; \bar u , \mu) $
for $\bar u\in H^1(\R^N)$ and $\mu \in [0,1]$.

\begin{proposition}\label{01092013-1435}
Suppose that {\em (\ref{18042014-0941})},  {\em (\ref{18042014-0949})} and {\em
(\ref{18042014-0950})}  are satisfied.\\
{\em (i)} For any bounded $U\subset H^1(\R^N)$ and $t>0$,
$\beta_{L^2} ({\bold \Psi}_t ( U \times [0,1])) \leq e^{- \bar v_{\infty} t} \beta_{L^2}(U)$;\\
{\em (ii)} If a bounded $U \subset H^1(\R^N)$ is relatively compact
as a subset of $L^2(\R^N)$, then
${\bold \Psi}_t (U\times [0,1])$ is relatively compact in $H^1(\R^N)$;\\
{\em (iii)} If $U \subset \overline{\conv}^{H^1} {\bold \Psi}_t(U \times
[0,1])$ for some bounded $U\subset H^{1}(\R^N)$ and $t>0$, then $U$ is relatively compact in $H^1(\R^N)$.
\end{proposition}
\noindent {\bf Proof:} It goes exactly along the lines of  \cite[Prop. 4.5]{Cwiszewski-RL}.

\section{Resonant averaging principle}

Let ${\bold A}: D({\bold A}) \to L^2(\R^N)$ be as in the previous section with ${\bold A_0}:=-\Delta$ and a continuous mapping ${\bold F}: [0, +\infty) \times H^1(\R^N) \to L^2(\R^N)$ satisfies conditions
(\ref{23052014-0558}) and (\ref{23052014-0559}).
\begin{theorem}\label{06012014-1625}
Assume that ${\cal N}:= \mathrm{Ker} {\bold A} \neq \{0\}$ and let $\bar{\bold F}: {\cal N} \rightarrow {\cal N}$ be given by
\begin{equation}\label{16032014-1420}
\bar{{\bold F}}(u):=\frac{1}{T}\int_0^T {\bold P} {\bold F} (s,u)\, \d s, \quad u \in {\cal N},
\end{equation}
where ${\bold P}: L^2(\R^N) \to N$ is the orthogonal projection onto ${\cal N}$ and
let open bounded sets $U \subset {\cal N}$ and $W \subset {\cal N}^\perp$ be such that $0 \notin \bar{{\bold F}}(\part U)$ and $0 \in W.$
By ${\bold \Phi}^{(\epsilon)}_T: H^1(\R^N) \rightarrow H^1(\R^N)$ denote the
translation operator by time $T$ for the problem
\begin{equation}\label{11032014-2223}
\dot{u}(t)= -{\bold A} u(t)+ \epsilon {\bold F} (t,u), \ t>0.
\end{equation}
Then there exists $\epsilon_0>0$ such that, for $\epsilon \in (0,\epsilon_0]$,
\begin{equation}
\Ind ({\bold \Phi}^{(\epsilon)}_T, U \oplus W)= (-1)^{m_-(\infty)}
\Deg_B (\bar{{\bold F}}, U), \nonumber
\end{equation}
where $m_-(\infty)$ is the total multiplicity of the negative eigenvalues of ${\bold A}$ and $\Deg_B$ stands for Brouwer topological degree.
\end{theorem}
\begin{remark} \label{22082014-1040}  {\em
(i) Let $\{e^{-{\bold A}t}: L^2(\R^N) \rightarrow L^2(\R^N)\}$ be the $C_0$-semigroup of bounded linear operators generated by
$-{\bold A}$. Then, for any $t>0$, ${\cal N}=\mathrm{Ker}\, {\bold A} = \mathrm{Ker}\,
(e^{-{\bold A}t}-I)$ (see \cite[Thm 16.7.2]{Hille}).\\
\indent (ii) Recall the known arguments on the spectrum of ${\bold A}$. Since the operator $-{\bold A}$ generates an analytic $C_0$ semigroup, the spectrum
$\sigma({\bold A})$ is contained in an interval $(-c,+\infty)$ with some $c>0$. Clearly $\sigma ({\bold A}_0 + {\bold V}_\infty) \subset [\bar v_\infty, +\infty)$. Since ${\bold V}_0({\bold A}_0+{\bold V}_\infty)^{-1}:
L^2(\R^N) \to L^2(\R^N)$ is a compact linear operator (see
\cite[Lem. 3.1]{Prizzi-FM}), by use of the Weyl theorem on essential spectra, we obtain $\sigma_{ess}({\bold A}) = \sigma_{ess}
({\bold A}_0-{\bold V}_0+{\bold V}_\infty)=\sigma_{ess}
({\bold A}_0+{\bold V}_\infty) \subset \sigma ({\bold A}_0 +{\bold V}_\infty) \subset [\bar v_\infty, +\infty)$. Hence, by general characterizations of essential
spectrum, we see that $\sigma({\bold A})\cap (-\infty,\bar v_\infty)$ consists of isolated eigenvalues with finite dimensional eigenspaces (see \cite{Schechter}).} In particular, $\dim {\cal N} < +\infty$.
\end{remark}
\noindent In the proof we shall need an auxiliary fact.
\begin{lemma} \label{07052014-1028}
Suppose that a sequence $(w_n(t))$ converges to the zero function in $L^2(\R^N)$ uniformly with respect to $t$ from compact subsets of $(0,T)$ and that the sequence $(v_n)$ of $T$-periodic solutions of the
problem
\begin{equation}
\dot{v}(t) = -{\bold A} v(t) + w_n(t), \qquad t \in [0,T] \nonumber
\end{equation}
is such that $v_n(0) \to \bar v_0$ in $L^2(\R^N)$ as $n \to +\infty$
for some $\bar v_0 \in H^1(\R^N)$. Then $v_n(t) \to v_0(t)$ in
$H^1(\R^N)$, uniformly with respect to $t$ from compact subsets of
$(0,T)$, where $v_0$ is a solution of
\begin{equation}
\left\{   \begin{array}{l}
\dot{v}(t)= -{\bold A} v(t),\qquad t \in [0,T]\\
v(0)= \bar v_0.
\end{array}  \right.
\end{equation}
\end{lemma}
\noindent {\bf Proof}. By the Duhamel formula it follows that, for $t \in [0,T]$,
\begin{equation}
v_n(t)-v_0(t) = e^{-{\bold A}t}(v_n(0)- \bar v_0) + \int_0^t e^{-{\bold A}(t-s)}\, w_n(s)\,\d s.\nonumber
\end{equation}
This implies that
\begin{equation}
\|v_n(t)-v_0(t)\|_{H^1} \leq C_{1/2} t^{-1/2}\|v_n(0)- \bar
v_0\|_{L^2} + C_{1/2} \int_0^t (t-s)^{-1/2}\|w_n(s)\|_{L^2} \d s
\nonumber
\end{equation}
for some constant $C_{1/2}>0$.
Let us take an arbitrary $\delta \in (0, T/2)$. Observe that
\begin{align}
\int_0^t \!\!\!(t-s)^{-1/2}\!\|w_n(s)\|_{L^2}\, \d s &\leq
\bigg(\frac{2}{\delta}\bigg)^{1/2}\!\!\!
\int_0^{t-\delta/2}\!\!\!\|w_n(s)\|_{L^2}\, \d s\! + \!
\int_{t-\delta/2}^t \!\!\!(t-s)^{-1/2}\|w_n(s)\|_{L^2}\, \d s \nonumber \\
& \leq \bigg(\frac{2}{\delta}\bigg)^{1/2}\!\!\! \int_0^{T-\delta/2}
\!\|w_n(s)\|_{L^2}\d s + 2 \bigg(\frac{2}{\delta}\bigg)^{-1/2}\!\!\!\!\!\! \sup_{s \in [\delta/2, T-\delta]}\|w_n(s)\|_{L^2}. \nonumber
\end{align}
Since $\|w_n(t)\|_{L^2} \to 0$ as $n \to \infty$ uniformly with
respect to $t$ from compact subsets of $(0,T)$ we infer that
$\|v_n(t) - v_0(t)\|_{H^1} \to 0$, as $n\to+\infty$, uniformly with respect to $t \in [\delta, T-\delta]$. \hfill $\square$

\noindent {\bf Proof of Theorem \ref{06012014-1625}}.
Let ${\bold \Theta}^{(\epsilon)}_T: H^1(\R^N) \times [0,1] \rightarrow
H^1(\R^N)$, $\epsilon \in [0,1]$, be the translation along trajectories operator for the problem
\begin{equation} \label{10032014-1000}
\dot{u}(t)= -{\bold A} u(t)+ \epsilon {\bold G} (t,u,\mu),  \qquad t>0,
\end{equation}
where ${\bold G}: [0,+\infty) \times H^1(\R^N) \times [0,1] \rightarrow L^2(\R^N)$ is the mapping given by the formula
\begin{equation}
{\bold G} (t,u,\mu):=(1-\mu){\bold F}(t,(1-\mu)u+\mu \tilde{{\bold P}}u)+\frac{\mu}{T}
\int_0^T {\bold P}{\bold F}(s,(1-\mu)u+\mu \tilde{{\bold P}}u)ds\nonumber
\end{equation}
with $\tilde{{\bold P}}: H^1(\R^N) \rightarrow {\cal N}$ being the orthogonal projection onto ${\cal N}$. Observe that, for $t \in [0,T]$ and $u \in
H^1(\R^N)$,
\begin{displaymath}
{\bold G}(t,u,0)= {\bold F}(t,u) \qquad \textnormal{and\ } \qquad {\bold G} (t,u,1)=
\bar{{\bold F}}(\tilde{{\bold P}}u).
\end{displaymath}
Clearly, note that, by (\ref{23052014-0558}),
for any $t, s \geq 0$, $u, v \in H^1(\R^N)$
and $\mu \in [0,1]$ one has
\begin{align}
\ & \| {\bold G} (t,u, \mu) - {\bold G} (s, v,\mu)\|_{L^2}  \leq C (1+ \|(1-\mu)u + \mu \tilde {\bold P} u \|_{H^1}) |t-s|^\theta + \nonumber \\
& \ \ \ \  + C \|(1-\mu)(u-v)+\mu \tilde {\bold P} (u-v)\|_{H^1} + \mu C\|(1-\mu)(u-v)+\mu \tilde {\bold P} (u-v)\|_{H^1} \nonumber \\
& \ \ \ \ \leq \tilde C (1+\|u\|_{H^1})|t-s|^\theta + \tilde C \|u-v\|_{H^1} \label{30042014-1209}
\end{align}
for some constants $\tilde C>0$ and $\theta \in (0,1)$. Further, in view of (\ref{23052014-0559}) and Remark (3.2)(b),
for any $t \geq 0$, $u \in H^1(\R^N)$ and $\mu \in [0,1]$,
\begin{align}
|{\bold G}(t, u, \mu)(x)| \!\leq\!  L(x)\left|[((1\!-\!\mu)I\!+\! \mu \tilde {\bold P})u](x)\right|\! +\! (K(x)+C K_0(x))(1\!+\!\|(1\!-\!\mu)u\!+\!\mu \tilde {\bold P}u\|_{H^1}) \ \label{30042014-1220}
\end{align}
where $K_0(x):=\sum_{k=1}^{\dim {\cal N}} |\varphi_k(x)|$ and $\{\varphi_k\}_{k=1}^{\dim {\cal N}}$ is an orthonormal basis of ${\cal N}$
with respect to the inner product induced in $L^2(\R^N)$ and some $C>0$.
Moreover, one immediately obtains that, for any $t \geq 0$, $u \in H^1(\R^N)$, $\mu,
\nu \in [0,1]$, $\left\| {\bold G}(t,u,\mu)-{\bold G}(t,u,\nu)\right\|_{L^2} \leq |\rho(\mu) - \rho(\nu)|(1+\|u\|_{H^1})$ for some $\rho\in C([0,1])$. Therefore ${\bold G}$ satisfies (\ref{18042014-0941}), (\ref{18042014-0949}) and (\ref{18042014-0950}). Hence it follows that ${\bold \Theta}^{(\epsilon)}_T$ is well defined and we can apply Proposition \ref{01092013-1435} to infer that ${\bold \Theta}^{(\epsilon)}_T$ is an ultimately compact operator (for any $\epsilon \in [0,1])$.\\
\indent Now we claim that there is $\epsilon_0 > 0$ such that, for any $\epsilon \in (0,\epsilon_0]$,
\begin{equation} \label{07032014-2307}
{\bold \Theta}^{(\epsilon)}_T (\bar u, \mu) \neq \bar u \qquad
\textnormal {for \ } \bar u \in \part(U \oplus W), \ \mu \in
[0,1].
\end{equation}
Suppose to the contrary that there are sequences $(\epsilon_n)$ in
$(0,+\infty)$, $(\bar u_n)$ in $\part(U \oplus W)$ and $(\mu_n)$ in $[0,1]$ such that $\epsilon_n \to 0^+$ and
\begin{equation} \label{04032014-2326}
{\bold \Theta}^{(\epsilon_n)}_T (\bar u_n, \mu_n) = \bar u_n \ \mbox{ for each } n \geq 1.
\end{equation}
Let $u_n: [0,+\infty) \to H^1(\R^N)$, $n \geq 1$ be solutions of (\ref{10032014-1000}) with $\epsilon =\epsilon_n$ and $\mu = \mu_n$, satisfying the initial value condition $u_n(0)=\bar u_n$.
By (\ref{04032014-2326}) and the Duhamel formula,
\begin{equation}
\bar u_n = u_n(T) = e^{-T{\bold A}} \bar u_n +
\epsilon_n \int_0^T e^{- (T-s){\bold A}}\, {\bold G} (s, u_n(s),\mu_n) \, \d s. \nonumber
\end{equation}
Moreover, observe that, for any $\phi \in {\cal N}$,
we have
\begin{align}
(\bar u_n, \phi)_{L^2} & = (e^{-T{\bold A}} \bar u_n,
\phi)_{L^2} + \epsilon_n \int_0^T
(e^{-(T-s){\bold A}}\, {\bold G}(s, u_n(s),\mu),
\phi)_{L^2}\, \d s  \nonumber \\
& = ({\bold P} e^{-T{\bold A}} \bar u_n,
\phi)_{L^2} + \epsilon_n \int_0^T ({\bold P} e^{- (T-s){\bold A}}\, {\bold G}(s, u_n(s),\mu),\phi)_{L^2}\, \d s. \nonumber
\end{align}
This, due to Remark (\ref{22082014-1040}) (i), yields
\begin{equation}
\label{07032014-1919} \int_0^T ({\bold G}(s, u_n(s),\mu_n), \phi)_{L^2} \, \d s = 0 \ \mbox{ for all } \phi\in {\cal N}.
\end{equation}
Furthermore, without loss of generality, we may assume that $\mu_n \to \mu_0$ as $n \to +\infty$ for some $\mu_0 \in
[0,1]$. By $(\ref{30042014-1220})$ it follows that there is $R>0$ such that $\|u_n(t)\|_{H^1} \leq R$ for all $t >0$ and $n \geq 1$. In view of Lemma \ref{08052014-1217}, for all  $m, n\geq 1$,
\begin{align}
\|(1-\chi_n)\bar u_n \|_{L^2}^2 &= \|(1-\chi_n)u_n(0) \|_{L^2}^2 =
\|(1-\chi_n) u_n (m T)\|_{L^2}^2 \nonumber \\ &\leq R^2
e^{-2\bar v_\infty m T} + \alpha_n
\end{align}
where $\chi_n$ is the characteristic function of $B(0,n)$. Since $m$ was arbitrary one gets $\|(1-\chi_n) \bar u_n\|_{L^2}\leq \sqrt{\alpha_n}$ for all $n\geq 1$. Due to the Rellich-Kondrachov, the sequence
$\left(\chi_n \bar u_n\right)$ is relatively compact in $L^2(\R^N)$. Therefore we can infer that $(\bar u_n )_{n\geq 1}$ is relatively compact in
$L^2(\R^N)$. And since it is bounded in $H^1(\R^N)$ we get a subsequence $(\bar u_{n_k})$ of $(\bar u_n)$, such that $\bar
u_{n_k} \to \bar u_0$ in $L^2(\R^N)$ for some $\bar u_0\in H^1(\R^N)$.
In view of Lemma \ref{07052014-1028}, $(u_n(t))$ converges in $H^1(\R^N)$ to $u_0(t)$ uniformly with respect to $t\in (0,+\infty)$ where $u_0$ is the solution of the problem
\begin{equation}
\dot{u}(t) = - {\bold A} u(t), \ t>0, \ \ u(0)=u(T)=\bar u_0. \nonumber
\end{equation}
By Remark \ref{22082014-1040} (i) we get $\bar u_0 \in {\cal N}$ and $u_0(t) =\bar u_0$ for all $t\in [0,+\infty)$.
Finally, after passing to the limit in (\ref{07032014-1919}), we
obtain that $\bar F(\bar u_0) =0$, a contradiction proving (\ref{07032014-2307}).\\
\indent Using (\ref{07032014-2307}) and the
homotopy invariance of the fixed point index for ultimately compact maps, we infer that, for  all $\epsilon\in  (0,\epsilon_0]$,
\begin{equation}
\Ind ({\bold \Theta}^{(\epsilon)}_T (\cdot, 0), U \oplus W ) = \Ind
({\bold \Theta}^{(\epsilon)}_T (\cdot, 1), U \oplus W). \nonumber
\end{equation}
Clearly
\begin{equation}
{\bold \Theta}^{(\epsilon)}_T (u,1) = {\bold \Psi}_T^{(\epsilon)}(\tilde {\bold P}u) + {\bold \Psi}_T(
(I-\tilde {\bold P})u). \nonumber
\end{equation}
where ${\bold \Psi}^{(\epsilon)}_T: {\cal N} \to {\cal N}$ is the translation along
trajectories operator associated with the equation
$\dot{u}(t)= \epsilon \bar {\bold F} (u(t))$, $t>0$, and the operator ${\bold \Psi}_T : {\cal N}^\perp \to {\cal N}^\perp$, where ${\cal N}^\perp$ stands for the space orthogonal to ${\cal N}$ in $H^1(\R^N)$, is given by
${\bold \Psi}_T (u):= e^{-T {\bold A}}u, \, u \in {\cal N}^\perp$. This means that ${\bold \Theta}^{(\epsilon)}_T (\cdot,1)$ is topologically adjoint with the mapping $\tilde {\bold \Theta}_{T}^{(\epsilon)}: {\cal N} \times {\cal N}^\perp \to {\cal N} \times {\cal N}^\perp$,
$\tilde {\bold \Theta}_{T}^{(\epsilon)} (u, v) := ({\bold \Psi}_T^{(\epsilon)} (u),
{\bold \Psi}_T(v))$, $u\in {\cal N}, v\in {\cal N}^{\perp}$. Therefore, by the product formula for fixed point index we get
\begin{equation}\label{02052014-14:52}
\Ind({\bold \Theta}_T^{(\epsilon)}(\cdot, 1), U \oplus W) = \Ind ({\bold \Psi}_T^{(\epsilon)}, U) \cdot \Ind ({\bold \Psi}_{T}, W).
\end{equation}
By the Krasnoselskii result (\cite[Lemma 13.1]{Krasnosielski}), decreasing $\epsilon_0$ if necessary, we get
\begin{equation}\label{02052014-14:55}
\Ind ({\bold \Psi}_T^{(\epsilon)} , U) = \Deg_B (\bar{\bold F}, U) \mbox{ for } \epsilon\in (0,\epsilon_0].
\end{equation}
To conclude we need to determine the fixed point index of $\Ind ({\bold \Psi}_{T}, W)$. According to Remark
(\ref{22082014-1040}) (ii), the set $\sigma({\bold A}) \cap (-\infty, 0)$ is bounded and closed. Consider the restriction $\tilde {\bold A}$ of ${\bold A}$ in the space $\tilde X$ orthogonal to ${\cal N}$ in $L^2(\R^N)$. Then, due to spectral theory (see \cite[Ch.7]{Dunford}),
there are closed subspaces $X_-$ and $X^+$ of $\tilde X$, such that $X_- \oplus X_+ = \tilde X$, $\tilde {\bold A} (X_-) \subset X_-$, $\tilde {\bold A} (D(\tilde {\bold A})\cap X_+) \subset X_+$, $\sigma (\tilde {\bold A}|_{X_-})=
\sigma({\bold A})\cap (-\infty,0)$, $\sigma(\tilde {\bold A}|_{X_+})= \sigma({\bold A}) \cap (0, +\infty)$. Define $ {\bold \Gamma}:
{\cal N}^{\perp} \times [0,1] \to {\cal N}^{\perp}$ by
$$
{\bold \Gamma} (\bar u, \mu):= e^{-T  {\bold A}} ((1-\mu) \bar u + \mu {\bold P}_- \bar u)
$$
where ${\bold P}_-: {\cal N}^{\perp} \to X_-$ is the restriction of the projection onto $X_-$ in $L^2(\R^N)$. Since $\dim
X_-<+\infty$ we infer that ${\bold P}_-$ is continuous. We also claim that
${\bold \Gamma}$ is ultimately compact. To see this take a bounded set $B
\subset H^1(\R^N)$ such that $B = \overline{\conv}^{H^1} {\Gamma} (
B\times [0,1]).$ This means that $B \subset \overline{\conv}^{H^1}
e^{-T  {\bold A}}(B \cup {\bold P}_- B)$. Since $B \cup {\bold P}_- B$ is bounded,
Proposition \ref{01092013-1435} (ii) implies that $B$ is relatively
compact in $H^1(\R^N)$, which proves the ultimate compactness of
${\bold \Gamma}$. Therefore, since $\mathrm{Ker} (I - {\bold \Gamma} (\cdot,\mu))=\{ 0\}$ for
$\mu\in [0,1]$, by the homotopy invariance and the restriction
property of the Leray-Schauder fixed point index, one gets
\begin{eqnarray}
\Ind ({\bold \Psi}_T, W ) & = & \Ind_{LS} ( e^{-T {\bold A}}{\bold P}_- , W ) \nonumber \\
& = & \Ind_{LS} ( e^{-T {\bold A|}_{X_-}}, W \cap X_-) = (-1)^{m_- (\infty)}.\label{02052014-1531}
\end{eqnarray}
The latter equality comes from the fact that $\sigma
({\bold A} |_{X_-})\subset (-\infty,0) $ consists of eigenvalues of finite dimensional eigenspaces. Finally, the proof is completed in view of (\ref{02052014-14:52}),
(\ref{02052014-14:55}) and (\ref{02052014-1531}). \hfill $\square$\\

\indent Using the above result and the existence property of fixed point index, one immediately obtains the following
existence result.
\begin{corollary}
Under the assumptions of Theorem {\em\ref{06012014-1625}}, if
$\Deg_B(\bar {\bold F}, U) \neq 0$, then for sufficiently small $\epsilon >0$, the problem {\em (\ref{11032014-2223})} admits a $T$-periodic solution.
\end{corollary}
\noindent We also derive the following {\em continuation principle}, which will be used in the proof of the main result.
\begin{theorem} \label{12032014-1608}
Under the assumptions of Theorem {\em \ref{06012014-1625}}, if for some $R_0 >0$ the following conditions are satisfied:
\begin{enumerate}
\item [{\em(i)}] $\Deg_B (\bar {\bold F}, B_{\cal N} (0, R_0)) \neq 0$;
\item [{\em (ii)}] for any $\epsilon \in (0,1)$ the problem {\em (\ref{11032014-2223})} has not a $T$- periodic solution with $\| u(0)\|_{H^1} \geq R_0$;
\end{enumerate}
then the equation
\begin{equation}
\dot{u}(t)= -{\bold A} u(t) + {\bold F}(t, u(t)), \qquad t>0 \nonumber
\end{equation}
admits a $T$-periodic solution.
\end{theorem}
\noindent {\bf Proof.} Take $U:= B_{\cal N} (0, R_0)$ and $W:= B_{{\cal N}^\perp}(0, R_0)$. By
Theorem \ref{06012014-1625}, there exists $\epsilon_0 >0$ such that for $\epsilon \in (0,
\epsilon_0]$
\begin{equation}
\Ind ({\bold \Phi}_T^{(\epsilon)}, U\oplus W)= (-1)^{m_-(\infty)} \Deg_B (\bar {\bold F},
U). \nonumber
\end{equation}
Note that $\part (U\oplus W) \subset H^1(\R^N) \setminus B_{H^1}(0, R_0)$, therefore by (ii), we deduce that, for any $\epsilon \in (0, 1)$ and $\bar u \in \part (U\oplus W)$, ${\bold \Phi}_T^{(\epsilon)}(\bar
u) \neq \bar u$. Hence, either ${\bold \Phi}_T^{(1)}$ has a fixed point in $\part (U\oplus W)$ (which proves the assertion) or, by the homotopy invariance of the fixed point index,
\begin{equation}
\Ind ({\bold \Phi}_T^{(1)},U\oplus W) = \Ind ({\bold \Phi}_T^{(\epsilon_0)}, U\oplus W) =
(-1)^{m_-(\infty)}\Deg_B (\bar {\bold F}, U)\neq 0, \nonumber
\end{equation}
which entails the existence of a fixed point in $U\oplus W$.\hfill $\square$

\section{Landesman-Lazer type criterion}

In this section we prove the main result -- Theorem \ref{15032014-2356}.  The proof is based on the continuation principle stated in Theorem \ref{12032014-1608}. Throughout this section we assume that $V$ and $f$ are as in Theorem \ref{15032014-2356}. Then the Nemytskii operator ${\bold F}$ determined by $f$ satisfies (\ref{23052014-0558}) and (\ref{23052014-0559}). The Landesman-Lazer type conditions stated in Theorem \ref{15032014-2356} imply the following inequalities.
\begin{lemma}\label{15032014-1047}
$\mbox{ }$
\begin{enumerate}
\item [\  {\em (i)}] If {\em (\ref{12032014-1611})} holds, then there exists $R_0 > 0$ such that
\begin{equation}\label{16032014-1500}
(\bar{{\bold F}}(u),u)_{L^2}>0 \textnormal{ for any } u \in {\cal N} \setminus
B_{\cal N} (0, R_0).
\end{equation}
\item[\  {\em (ii)}]
If {\em (\ref{12032014-1618})} holds, then there exists $R_0 >0$ such that
\begin{equation}\label{16032014-1452}
(\bar{\bold F}(u),u)_{L^2}< 0 \textnormal { for any } u \in {\cal N} \setminus
B_{\cal N} (0, R_0).
\end{equation}
\end{enumerate}
\end{lemma}
\noindent {\bf Proof}. (i) Suppose to the contrary that there exists a sequence
$(\bar u_n)$ in ${\cal N}$ such that $\|\bar u_n\|_{H^1} \to +\infty$ as $n \to +\infty$ and $(\bar {\bold F}(\bar u_n),\bar u_n)_{L^2} \leq 0$. Put
$\mu_n := \|\bar u_n \|_{H^1}$ and let $\bar v_n:= \bar u_n/ \mu_n$. Clearly, $(\bar v_n)$ is bounded in $H^1(\R^N)$ and, since $\dim {\cal N} < +\infty$, we may assume that, up to a subsequence, $\bar v_n \to \bar v_0$ as $n \to +\infty$ in $H^1(\R^N)$ for some $\bar v_0 \in {\cal N}$, which implies $\bar v_n \to \bar
v_0$ in $L^2(\R^N)$. On the other hand
\begin{align}
0 &\geq (\bar {\bold F} (\bar u_n), \bar u_n)_{L^2} = \frac{1}{T}
\int_0^T \left( {\bold P} {\bold F} (t, \bar u_n), \bar u_n \right)_{L^2}\d t \nonumber \\
& = \frac{1}{T} \int_0^T ({\bold F}(t, \bar u_n), \bar u_n)_{L^2}\,\d t =
\frac{\mu_n}{T} \int_0^T \int_{\R^N} f (t,x,\mu_n v_n(x)) v_n(x)\, \d x
\d t. \nonumber
\end{align}
Again passing to subsequence, we assume that $\bar
v_n(x) \to \bar v_0(x)$ for almost every $x \in \R^N$ and that there is $k \in L^1(\R^N)$ such that, for all $n\geq 1$, $\left|\bar v_n(x) \right| \leq k(x)$ for almost every $x \in \R^N$. In view of the Fatou lemma, we obtain
\begin{align}
0 & \geq \liminf_{n \to+\infty} \int_0^T \int_{\R^N} f
(t,x,\mu_n \bar v_n(x)) \bar v_n(x)\, \d x \d t \nonumber \\
& \geq \int_0^T \bigg(\int_{\{\bar v_0 >0\}}\check f_+(t,x) \bar
v_0(x)\, \d x+ \int_{\{\bar v_0 <0\}} \hat f_-(t,x) \bar v_0(x)\, \d x
\bigg)\, \d t>0 ,\nonumber
\end{align}
a contradiction with the assumption (\ref{12032014-1611}) that proves the assertion (i).\\
\indent (ii) can be proved analogically. \hfill $\square$\\

\noindent {\bf Proof of Theorem \ref{15032014-2356}}.
We claim that if  $R_0>0$ is as in Lemma \ref{15032014-1047}  then, for all $R>R_0$,
$$
\Deg_B (\bar {\bold F}, B_{\cal N}(0,R)) = \left\{\begin{array}{cl}
1 & \mbox{ if (\ref{12032014-1611}) holds},\\
(-1)^{\dim {\cal N}} &  \mbox{ if  (\ref{12032014-1618}) holds}.
\end{array} \right.
$$
Indeed, assume first (\ref{12032014-1611}) and define ${\bold H}:D_{\cal N} (0,R)\times [0,1]\to {\cal N}$, ${\bold H}(u,\mu):=\mu \bar {\bold F}(u)+(1-\mu)u$,
$u\in D_{\cal N}(0,R)$, $\mu\in[0,1]$. It is clear that
${\bold H}(\cdot, \mu)$ has no zeros in the boundary
$\part D_{\cal N}(0,R)$ due to (\ref{16032014-1500}), which, by use of the
homotopy invariance property, yields
$\Deg_B ({\bar{\bold F}}, B_{\cal N}(0,R))=1$. In a similar manner we show that if (\ref{12032014-1618}), then
the mapping ${\bold H}:D_{\cal N} (0,R)\times [0,1]\to {\cal N}$, ${\bold H}(u,\mu):=\mu \bar {\bold F}(u)-(1-\mu)u$,
$u\in D_{\cal N}(0,R)$, $\mu\in[0,1]$, has no zeros in
$\part D_{\cal N}(0,R)$, which gives
$\Deg_B(\bar{\bold F}, B_{\cal N}(0,R))=(-1)^{\dim {\cal N}}$.\\
\indent Now we claim that there is $R_0>0$ such that the problem
\begin{equation} \label{16032014-0008}
\dot{u}(t) = -{\bold A} u(t) + \epsilon {\bold F}(t,u),
\, t>0,
\end{equation}
has no $T$-periodic solutions for $\epsilon \in (0,1)$ with $\|u(0)\|_{H^1} \geq R_0$. Suppose to the contrary that there  are $\epsilon_n\in (0,1)$ and $T$-periodic solutions $u_n:[0,T] \to H^1(\R^N)$ of (\ref{16032014-0008}) with $\epsilon :=\epsilon_n$, $n\geq 1$, such that $\|u_n(0)\|_{H^1} \to +\infty$ as $n \to +\infty$. Put $\mu_n:=\sup_{t\geq 0} \|u_n(t)\|_{H^1}$ and let $v_n:= \mu_n^{-1} u_n$. Then one can easily observe that $v_n$ is a $T$-periodic solution of
\begin{equation}
\dot{v}(t) = - {\bold A} v(t) + {\bold F}_n (t,v(t)), \ \ t \in [0,T],
\label{02052014-1102}
\end{equation}
with ${\bold F}_n(t,u):=\epsilon_n \mu_n^{-1} {\bold F} (t,\mu_n u)$, $t\geq 0$, $u\in H^1(\R^N)$.
Clearly, by use of (\ref{23052014-0558}) and (\ref{23052014-0559}), for sufficiently large $n$ and all $t,s\in [0,1]$, $u,v\in H^1(\R^N)$, we have
\begin{align}
\|{\bold F}_n (t,u) -{\bold F}_n (s,v)\|_{L^2} \leq \epsilon_n \mu_n^{-1} C(1+\|\mu_n u\|_{H^1})|t-s|^{\theta} + \epsilon_n \mu_n^{-1}C\|\mu_n u-\mu_n v\|_{H^1} \nonumber\\
\leq C(1+\|u\|_{H^1})|t-s|^\theta+C\|u-v\|_{H^1}, \nonumber
\end{align}
\begin{align}
| {\bold F}_n (t,u)(x) |\leq \epsilon_n \mu_n^{-1}
L(x)|\mu_n u(x)| + \epsilon_n \mu_n^{-1}K(x)(1+\|\mu_n u\|_{H^1}) \nonumber \\
\leq L(x)|u(x)|+K(x)(1+\|u\|_{H^1}) \ \ \mbox{ for a.a. } x\in \R^N. \nonumber
\end{align}
Hence, by Lemma \ref{08052014-1217}, for all $m,n\geq 1$,
$$
\|(1-\chi_n) v_n(0)  \|_{L^2}^2 = \|(1-\chi_n) v_n (m T)\|_{L^2}^2 \leq \tilde R^2 e^{-2\bar v_\infty  m T} +\alpha_n
$$
where $\chi_n$ is the characteristic function of $B(0,n)$, $\alpha_n \to 0^+$ as $n\to+\infty$ ($\alpha_n$ depends only on $V$ and $K$, $L$, which are common for all ${\bold F}_n$) and $\tilde R>0$ such that $\|v_n (t)\|_{H^1} \leq \tilde R$ for all $t\geq 0$ and $n\geq 1$. Since $m$ is arbitrary we see that $\|(1-\chi_n)  v_n(0)\|_{L^2}\leq \sqrt{\alpha_n}$ for $n\geq 1$. Due to the Rellich-Kondrachov,  $\{\chi_n v_n(0) \}_{n\geq 1}$ is relatively compact in $L^2(\R^N)$. Therefore $\{ v_n(0) \}_{n\geq 1}$ is relatively compact in $L^2(\R^N)$. As a bounded sequence in $H^1(\R^N)$ it contains a subsequence convergent in $L^2(\R^N)$ to some  $\bar v_0\in H^1(\R^N)$, therefore we assume that $v_n(0)\to \bar v_0$ in $L^2(\R^N)$.\\
\indent Moreover, for all $t\geq 0$,
$$
\|{\bold F}_n(t,v_n(t)) \|_{L^2}
\leq \|\chi_m {\bold F}_n (t,v_n(t))\|_{L^2}+
\|(1-\chi_m) {\bold F}_n (t,v_n(t))\|_{L^2}.
$$
Since $f$ is bounded, it is clear that ${\bold F}_n(t,v_n(t))(x) \to 0$ as $n\to+\infty$, for a.a. $x\in \R^N$, which gives, for each $m\geq 1$,
\begin{equation} \label{25082014-1605}
\max_{t\geq 0} \|\chi_m {\bold F}_n (t,v_n(t))\|_{L^2} \to 0 \mbox{ as } n\to +\infty.
\end{equation}
Furthermore, for all $m,n\geq 1$ and $t\geq 0$,
\begin{align}
\|(1\!-\!\chi_m) {\bold F}_n (t,v_n(t))\|_{L^2} \leq
C^{N/p}\|(1\!-\!\chi_m)L\|_{L^p} \|v_n(t)\|_{H^1} \!+\!
\|(1\!-\!\chi_m)K \|_{L^2} (1\!+\!\|v_n(t)\|_{H^1}) \nonumber\\
\leq \beta_m:= C^{N/p} \|(1\!-\!\chi_m)L\|_{L^p} \!+\!2\|(1\!-\!\chi_m)K \|_{L^2}
\nonumber
\end{align}
where $C>0$ is the constant related to the embedding
$H^1(\R^N) \subset L^{2N/(N-2)} (\R^N)$. This, together with (\ref{25082014-1605}), gives
$$
\limsup_{n\to+\infty} \max_{t\geq 0} \|{\bold F}_n (t,v_n(t))\|_{L^2} \leq \beta_m
$$
and since $\beta_m\to 0^+$ as $m\to +\infty$,
we get $\max_{t\geq 0} \|{\bold F}_n(t, v_n(t))\|_{L^2}\to 0$ as $n\to +\infty$. Hence, by Lemma \ref{07052014-1028}, we infer that
$(v_n)$ converges in $C([0,T],H^1(\R^N))$ to some $v_0:[0,T]\to H^1(\R^N)$ being the $T$-periodic solution of
\begin{equation}
\dot{v}(t) = -{\bold A}v(t), \ v(0)=\bar v_0. \nonumber
\end{equation}
This means that $\bar v_0 = e^{-{\bold A }T} \bar v_0$, i.e. $v_0 (t)= \bar v_0$ for $t\geq 0$, and, since
$\max_{t\geq 0} \|v_n(t)\|_{H^1}=1$ for any $n
\geq 1$, we have $\bar v_0 \neq 0$.\\
\indent On the other hand, by the $T$-periodicity of $v_n$ and by the Duhamel formula it follows that
\begin{equation}
(v_n(T), \bar v_0)_{L^2} = (e^{-{\bold A}T}v_n(0), \bar v_0)_{L^2} + \epsilon_n \mu_n^{-1} \int_0^T (e^{-{\bold A}(T-t)}\,{\bold F} (t, \mu_n v_n(t)), \bar
v_0)_{L^2}\, \d t, \nonumber
\end{equation}
and
\begin{equation}
(v_n(0), \bar v_0)_{L^2}= (v_n(0), e^{-{\bold A}T} \bar v_0)_{L^2}+\epsilon_n
\mu_n^{-1} \int_0^T ({\bold F}(t, \mu_n v_n(t) ), e^{-{\bold A}(T-t)} \bar v_0)_{L^2}\,
\d t,\nonumber
\end{equation}
i.e., for all $n\geq 1$,
\begin{equation}
\int_0^T ({\bold F}(s, \mu_n v_n(t)), \bar v_0)_{L^2} \d t = 0. \nonumber
\end{equation}
\indent Assume now that (\ref{12032014-1611}) holds. Then, by the Fatou lemma,
\begin{equation} \label{16032014-1259}
0 = \liminf_{n \to +\infty} \int_0^T ({\bold F}(t, \mu_n v_n(t)), \bar
v_0)_{L^2} \d t \geq \int_0^T \bigg( \liminf_{n \to +\infty}({\bold F}(t,
\mu_n v_n(t)), \bar v_0)_{L^2}\bigg) \d t.
\end{equation}
Fix any $t\in [0,T]$ and let $(n_k)$ be an increasing sequence of positive integers such that
\begin{equation}\label{16032014-1253}
\liminf_{n \to +\infty} ({\bold F}(t, \mu_n v_n(t)), \bar v_0)_{L^2} =
\lim_{k \to + \infty} ({\bold F}(t, \mu_{n_k}v_{n_k}(t)), \bar v_0)_{L^2}
\end{equation}
and $(v_{n_k}(t))$ converges to $\bar v_0$ almost everywhere (the set on which convergence occurs may depend on $t$). Again due to the Fatou lemma it follows that
\begin{align}\label{16032014-1252}
\lim_{k \to +\infty} ({\bold F} (t, \mu_{n_k} v_{n_k}(t)), \bar v_0)_{L^2}
&\geq \int_{\R^N} \liminf_{k \to +\infty}
f(t,x,\mu_{n_k}v_{n_k}(t)(x)) \bar v_0(x)\, \d x \nonumber \\
&\geq \int_{\{\bar v_0 >0\}}\check f_+(t,x) \bar v_0(x)\, \d x +
\int_{\{\bar v_0 < 0\}}\hat f_-(t,x) \bar v_0(x)\, \d x,
\end{align}
since for almost every $x \in \{\bar v_0 >0\}$,
\begin{equation}
\liminf_{k \to +\infty} f(t,x,\mu_{n_k} v_{n_k}(t)(x)) \geq \check
f_+ (t,x) \nonumber
\end{equation}
and, for almost every $x \in \{\bar v_0 < 0\}$,
\begin{equation}
\limsup_{k \to +\infty} f(t,x,\mu_{n_k} v_{n_k}(t)(x)) \leq \hat f_-
(t,x). \nonumber
\end{equation}
Summing up, by (\ref{16032014-1253}) and (\ref{16032014-1252}), we get, for any $t \in [0,T]$,
\begin{equation}
\liminf_{k\to+\infty}({\bold F}(t,\mu_nv_n(t)), \bar v_0)_{L^2} \geq
\int_{\{\bar v_0>0\}}\check f_+(t,x) \bar v_0(x)\,  \d x + \int_{\{\bar
v_0<0\}}\hat f_-(t,x) \bar v_0(x)\, \d x, \nonumber
\end{equation}
which together with (\ref{16032014-1259}) gives
\begin{align}
0 &= \liminf_{k\to+\infty}\int_0^T (F(t, \mu_n v_n(t)), \bar
v_0)_{L^2}\, \d t \nonumber \\ &\geq \int_0^T \int_{\{\bar
v_0>0\}}\check f_+(t,x) \bar v_0(x)\, \d x \d t + \int_0^T \int_{\{\bar
v_0<0\}}\hat f_-(t,x) \bar v_0(x)\, \d x \d t >0, \nonumber
\end{align}
a contradiction proving the assertion, i.e. the condition (ii) of  Theorem \ref{15032014-2356} in  the case (\ref{12032014-1611}) holds.\\
\indent The case when (\ref{12032014-1618}) is satisfied can be treated in a analogous manner.\\
\indent We complete the proof by using Theorem \ref{12032014-1608}. \hfill $\square$

\end{document}